\documentclass[letterpaper, 10pt, journal, twocolumn, final]{IEEEtran}

\usepackage{mathtools}
\usepackage{bm}
\usepackage{algorithmicx}
\usepackage{algorithm}
\usepackage{algpseudocode}
\usepackage{graphicx}
\usepackage{epstopdf}
\usepackage{amsfonts}
\usepackage{amsmath}%
\usepackage{dsfont}

\usepackage{epsfig} 
\usepackage{xcolor}
\usepackage{amssymb}

\usepackage{cuted}

\usepackage{enumitem}
\usepackage{tikz}
\usepackage{balance}

\usetikzlibrary{shapes}
\usetikzlibrary{positioning}

\usepackage[hidelinks]{hyperref}

\newtheorem{theorem}{Theorem}

\newtheorem{lemma}{Lemma}

\newtheorem{algo}{Algorithm}

\newtheorem{assumption}{Assumption}

\newtheorem{remark}{Remark}

\newcommand\oprocendsymbol{\hbox{$\square$}}
\newcommand\oprocend{\relax\ifmmode\else\unskip\hfill\fi\oprocendsymbol
}

\usepackage[normalem]{ulem} %

\newcommand{\R}{{\mathbb{R}}}
\renewcommand{\natural}{{\mathbb{N}}}

\newcommand{\until}[1]{\{1,\ldots,#1\}} 

\newcommand{\1}{\mathbf{1}}

\newcommand{\norm}[1]{\left \|#1 \right \|}

\newcommand{\cD}{\mathcal{D}}

\newcommand{\real}{\mathbb{R}}

\newcommand{\av}{_{\text{avg}}}

\newcommand{\tzperp}{\tilde{z}_\perp}

\newcommand{\col}{\mathtt{col}}
\newcommand{\N}{\mathbb{N}}
\newcommand{\B}{\mathcal{B}}

\newcommand{\tzm}{\tilde{z}_{\text{m}}}
\newcommand{\per}{\tau_{\text{per}}}

\definecolor{blue@O4S}{RGB}{0, 41, 69}
\definecolor{emph@O4S}{RGB}{0, 93, 137}
\definecolor{red@O4S}{RGB}{127,0,0}
\definecolor{gray@O4S}{RGB}{112, 112, 112}

\usepackage[printonlyused]{acronym}

\usepackage{todonotes}

\usepackage{tikz}
\usetikzlibrary{shapes,arrows}
\usetikzlibrary{fit}
\usepackage{expl3}

\graphicspath{{figs/}}

\def\nAg{N}
\def\xD{n}

\def\algo/{Extremum Seeking Tracking}

\begin{document}

\tikzstyle{triangle} = [draw, regular polygon, isosceles triangle]

\title{Extremum Seeking Tracking\\ for Derivative-free
  Distributed Optimization}

\author{Nicola Mimmo, Guido Carnevale, Andrea Testa, and Giuseppe Notarstefano
 \thanks{This work was supported in part by the Italian Ministry of Foreign Affairs and International Cooperation”, grant number BR22GR01.
 The authors are with the Department of Electrical and Information Engineering,  "Guglielmo Marconi", University of Bologna, 40126 Bologna, Italy, e-mail: \{nicola.mimmo2, guido.carnevale, a.testa, giuseppe.notarstefano\}@unibo.it.
  }}

\maketitle

\begin{abstract}
  In this paper, we deal with a network of agents that want to
  cooperatively minimize the sum of local cost functions depending on
  a common decision variable. We consider the challenging scenario in
  which objective functions are unknown and agents have only access to
  local measurements of their local functions. We propose a novel
  distributed algorithm that combines a recent gradient tracking
  policy with an extremum seeking technique to estimate the global
  descent direction. The joint use of these two techniques results in
  a distributed optimization scheme that provides arbitrarily accurate
  solution estimates through the combination of Lyapunov and averaging
  analysis approaches with consensus theory. We perform numerical
  simulations in a personalized optimization framework to corroborate
  the theoretical results.
\end{abstract}

\section{Introduction}

In recent years, distributed optimization over networks has become a
key research topic, see, e.g.,
\cite{nedic2018distributed,yang2019survey,notarstefano2019distributed}
for an overview, also in the setting with partially unknown cost function~\cite{conn2009introduction}.  
Examples include data analytics in machine learning~\cite{liu2020primer} as well as automatic
controller tuning~\cite{zalkind2020automatic,maass2021zeroth}.
We organize the literature review in two blocks: \emph{collaborative/distributed extremum seeking} and \emph{distributed zeroth-order/derivative-free optimization}.
Early references on distributed extremum seeking for the so-called
\emph{consensus optimization} framework
are~\cite{Menon2014Collaborative} for a discrete-time setting
and~\cite{Ye2016Distributed} for a continuous-time one.
A proportional-integral extremum seeking design technique is proposed
in~\cite{guay2018distributed}, while, in~\cite{Dougherty2017Extremum},
the gradient is approximated through a real-time
protocol. In~\cite{Salamah2018Cooperative}, authors propose the use of
the sliding mode to generate the dither signal at the base of the
extremum seeking.
More recently, in~\cite{li2020cooperative}, a distributed stochastic extremum
seeking scheme is proposed for a source localization problem, while
in~\cite{guay2021distributed} a distributed scheme approximating a
Newton method is proposed.
As for constraint-coupled distributed optimization,
in~\cite{Poveda2013Distributed}, a distributed extremum seeking
control, based on evolutionary game theory, is designed for
real-time resource allocation.
In~\cite{poveda2017distributed}, instead, a distributed,
continuous-time scheme based on sign-based consensus is designed.
In~\cite{Michalowsky2017Distributed}, a Lie bracket technique and extremum seeking are used for problems with linear constraints. In~\cite{Wang2019Distributeda}, resource
allocation problems are addressed by an extremum seeking scheme.

As for \emph{distributed zeroth-order/derivative-free optimization}, 
authors in~\cite{mhanna2022zero} develop a zeroth-order scheme based on a
1-point estimator and a gradient tracking policy. The
work~\cite{tang2020distributed} instead proposes a zeroth-order algorithm based on a
two-point estimator with a distributed gradient descent strategy and another one based on an $\xD$-point estimator with a gradient tracking policy.
Authors in~\cite{Lu2012Zero} propose a continuous-time gradient-free approach emulating
a distributed gradient algorithm for which optimal asymptotic convergence is
guaranteed. In~\cite{Liu2014Sample}, a sampled version of \cite{Lu2012Zero} is
proposed. In~\cite{Yuan2015Randomized}, a continuous-time distributed algorithm
based on random gradient-free oracles is proposed for convex optimization
problems. In~\cite{pang2019randomized}, the randomized gradient-free oracles
introduced in~\cite{Yuan2015Randomized} are used to build a gradient-free
distributed algorithm in directed networks.
Randomized gradient-free algorithms are used also in~\cite{Chen2017Randomized},
where sequential Gaussian smoothing is used for non-smooth distributed convex
constrained optimization. In~\cite{Ding2017Distributed}, gradient-free
optimization is addressed with the additional constraint that each agent can
only transmit quantized information.
The work in~\cite{Pang2018Exact} instead develops a ``directed-distributed
projected pseudo-gradient'' descent method for directed
graphs. Paper~\cite{Wang2018Design} combines the gradient-free strategy
of~\cite{Yuan2015Randomized} with a saddle-point algorithm. Authors
in~\cite{Wang2019Distributed} address an online constrained optimization problem
by relying on the Kiefer-Wolfowitz algorithm to approximate the gradients,
and~\cite{Bilenne2020Fast} combines the estimation of the gradient via a
``simultaneous perturbation stochastic approximation'' technique with the so-called matrix exponential learning optimization
method. In~\cite{Pang2020Randomized}, a randomized gradient-free method is
combined with a state-of-the-art distributed gradient descent approach for
directed networks. In~\cite{Sahu2020Decentralized}, an overview of zeroth-order
methods based on a Frank-Wolfe framework is provided. Authors
in~\cite{Wang2020Random} propose a distributed random gradient-free protocol to
solve constrained optimization problems by using projection techniques.
All the cited papers but~\cite{Bilenne2020Fast} prove stability of the
proposed schemes by using arguments not based on Lyapunov theory.
We also note that~\cite{Bilenne2020Fast} proposes an algorithm neither based on gradient tracking nor using extremum seeking.
This paper proposes a novel algorithm to solve a distributed optimization
problem in which network agents can only evaluate their local cost function at a
given point, but not its gradient. 
The proposed solution consists of a distributed protocol in which
	the (unavailable) local gradients are approximated through an extremum
seeking scheme.
The approximations of the gradients are used to feed a suitable tracking mechanism which, in turn, allows for steering the local solution estimates along approximations of the global descent direction.
The distributed algorithm uses a further consensus action on the local solution estimates.
The convergence of our scheme is proved through Lyapunov and averaging
theories for discrete-time systems.
It is worth mentioning that our scheme, together with the ones
in~\cite{Vandermeulen2018Discrete,kvaternik2012analytic}, is the only distributed extremum seeking scheme proposed in
discrete-time with the following distinctive features. 
The work~\cite{Vandermeulen2018Discrete} (i) does not address a consensus optimization problem and (ii) relies on consensus dynamics estimating the global cost, while ours estimates the global gradient.
Instead, in \cite{kvaternik2012analytic}, (i) the addressed consensus optimization problems have scalar decision variables and (ii) an extremum seeking technique is combined with a distributed gradient algorithm, i.e., without a tracking mechanism.
As for the literature on distributed zeroth-order/derivative-free methods, the closest work is~\cite{mhanna2022zero} in which extremum seeking is not used and the global gradient is approximated via a randomized 1-point policy.
Consistently with the comparison table provided in~\cite{mhanna2022zero}, we highlight that there are no other distributed algorithms in the literature using 1-point gradient estimators. Indeed, although the distributed scheme in~\cite{li2021distributed} uses 1-point gradient estimators, we remark that it is tailored for partition-based optimization, i.e., a simplified setup in which each local function depends only on the neighbors' decision variables, see~\cite{notarstefano2019distributed}.
Also, our Lyapunov-based tools combined with averaging theory for
discrete-time systems represent a distinctive feature in the algorithm
analysis.
The paper unfolds as follows. 
Section \ref{sec:ProblemFormulation} introduces
the problem and the proposed algorithm. 
The main result is provided in Section~\ref{sec:ProposedSolution} and proved in Section
\ref{sec:MainResult}.
In Section~\ref{sec:Simulation}, we test the proposed algorithm via numerical simulations.%

\subsubsection*{Notation}
Given $N$ vectors $x_1,\dots,x_N \in \mathbb{R}^n$, we denote by $\mathtt{col}(x_1,\dots,x_N)$
their column stacking. 
Given $N$ scalars $d_1, \dots, d_N$, we denote by $\mathtt{diag}(d_1, \dots, d_N)$ the diagonal matrix with $i$-th entry $d_i$.
The Kronecker product is denoted by $\otimes$. 
The identity matrix in $\R^{n\times n}$ is $I_n$. 
 The column vector of $N$ ones is denoted by $1_N$ and we define
$\1 := 1_N \otimes I_n$.  
Dimensions are omitted whenever they are clear
from the context. 
Finally, for $r >0$, we let $\B_r=\{x \in \R^{n} \,:\, \norm{x} \leq r\}$.

\section{Problem Formulation} %
\label{sec:ProblemFormulation}

We consider a network of $N$ agents communicating according to an undirected graph
$\mathcal{G}=(\mathcal{V},\mathcal{E},\mathcal{A})$, where $\mathcal{V} = \{1,\dots,N \}$
is the set of agents, $\mathcal{E} \subseteq \mathcal{V} \times \mathcal{V}$ is
the set of edges, and $\mathcal{A} =[a_{ij}] \in \R^{N \times N}$ is the weighted adjacency matrix. 
Agent $i$ and $j$ can exchange information only if $(i,j) \in \mathcal{E}$.
Accordingly, it holds $a_{ij} \ge 0$ if $(i,j) \in \mathcal{E}$ and $a_{ij} = 0$ otherwise.
We denote as $\mathcal{N}_i = \{j\in\mathcal{V}\mid (i,j)\in\mathcal{E}\}$ the set of neighbors of agent $i$.  
Moreover, we also associate to the graph the Laplacian matrix
  $\mathcal{L} := D - \mathcal{A} \in \R^{N \times N}$, where
  $D := \mathtt{diag}(\mathtt{deg}_1,\dots,\mathtt{deg}_N) \in \R^{N \times N}$ is the so-called degree matrix in which
  $\mathtt{deg}_i := \sum_{j \in \mathcal{N}_i}a_{ij}$ is the degree of agent $i$.
The next assumption specifies the class of considered graphs.
\begin{assumption}\label{ass:connectivity}
  The graph $\mathcal{G}$ is connected and the adjacency matrix $\mathcal{A}\in\real^{N\times N}$ is symmetric. \oprocend
\end{assumption}
In the proposed distributed setup, each agent $i$ is equipped with a sensor only
providing measurements of the local cost function $f_i: \R^\xD \to \R$ and aims at solving %
the problem%
\begin{align}
\label{eq:TotalCost}
\min_{w \in \mathbb{R}^\xD}\sum_{i=1}^{\nAg} f_i(w).
\end{align}
We enforce the following assumptions about the problem.%

\begin{assumption}\label{ass:strong_convexity}
  For all $i \in \until{N}$, the function $f_i$ is $\underline{L}$-strongly convex for some
  $\underline{L}>0$. \oprocend
\end{assumption}
\begin{assumption}\label{ass:lipschitz}
Each cost $f_i$ is $\mathcal{C}^3$ and 
has $\overline{L}_i$-Lipschitz continuous gradients.
We denote $\overline{L} := \max\{\overline{L}_1,\dots,\overline{L}_\nAg\}$. \oprocend
\end{assumption}
  Assumption~\ref{ass:strong_convexity} ensures the existence of a unique solution $x^\star\in\mathbb{R}^\xD$ to problem~\eqref{eq:TotalCost}.
Our aim is to iteratively find it via a distributed algorithm.
Namely, given the iteration index $t \in \N$ and by denoting with $w_i^t \in \R^{\xD}$ the $i$-th agent's estimate, at iteration $t$, of the solution to problem~\eqref{eq:TotalCost}, our goal is to design a distributed protocol able to steer all these estimates to the minimizer $x^\star$. 
The peculiar challenge of this paper is that each agent $i$ cannot access either the gradients $\nabla f_i(w_i^t)$ (as in standard gradient-based methods) or the cost functions in arbitrary points (as in standard zeroth-order methods).
More in detail, we assume that agent $i$ can only use the single measurement $f_i(w_i^t)$ properly combined with the so-called dither signal $d_i^t \in \R^{\xD}$ and the amplitude parameter $\delta > 0$.
The role of $d_i^t$ and $\delta$ will become clearer in the next section.

\section{\algo/:\\ Algorithm Introduction and Convergence}
\label{sec:ProposedSolution}

Being the gradients $\nabla f_i$ not available, we replace them with an estimation based on a proper elaboration of the local costs excited via suitable dithering signals $d_i^t$.
The arising distributed method, termed \algo/, is described in Algorithm~\ref{table:algorithm}, from agent $i$ perspective.
In Algorithm~\ref{table:algorithm}, $\gamma >0 $ represents the step size while the parameter $\delta > 0$ represents the amplitude of the
dither signal $d_i^t$
defined as
\begin{align}
  \label{eq:dither}
  d_i^t = \mathtt{col}\bigg(\sin\bigg(\dfrac{2\pi t}{\tau_{i_1}} + \phi_{i_1}\bigg),\, \dots, \, \sin\bigg(\dfrac{2\pi t}{\tau_{i_\xD}} + \phi_{i_\xD}\bigg)\bigg),
\end{align}
where $\tau_{i_p} \in \mathbb{N}$ %
and $\phi_{i_p} \in \mathbb{R}$ such that, given $p,q,r \in \until{\xD}$,
$p \ne q$, $q \ne r$, $p \ne r$, it holds
\begin{subequations}\label{eq:frequencies}
	\begin{align}
		&	\sum_{t=0}^{\per-1} \sin\left(\tfrac{2\pi t}{ \tau_{i_p}} + \phi_{i_p}\right) \! = \! 0
		\label{eq:frequencies1}
		\\
		&	\sum_{t=0}^{\per-1} \sin\left(\tfrac{2\pi t}{ \tau_{i_p}} + \phi_{i_p}\right) \sin\left(\tfrac{2\pi t}{\tau_{i_q}} + \phi_{i_q}\right) \! = \! \dfrac{\per}{2}
			\label{eq:frequencies2}
		\\
		&	\sum_{t=0}^{\per-1} \sin\!\left(\tfrac{2\pi t}{ \tau_{i_p}} \! + \! \phi_{i_p}\right)\! \sin\!\left(\tfrac{2\pi t}{\tau_{i_q}}\! + \! \phi_{i_q}\right)\!\sin\!\left(\tfrac{2\pi t}{ \tau_{i_r}} \! + \! \phi_{i_r}\right) \!=\! 0,\label{eq:frequencies3}
	\end{align}
\end{subequations}
for all $i \in \until{N}$. 
Here, $\per \in \mathbb{N}$ is the least common multiple of all
periods $\tau_{i_p}$.
Sinusoidal dither functions are useful in practical
applications to guarantee smooth inputs to the plant. However, other
possibilities, e.g., as square or triangular waves, are possibile~\cite{TAN2008choice}.
\begin{algorithm}[t]
	\begin{algorithmic}
		\State initialization: $x_i^0 \in \mathbb{R}^n$ and $z_{i}^0 = 0$
		\For{$t=0, 1, \dots$}
								\vspace{-3ex}
		\State
		\begin{subequations}\label{eq:GTA}
			\begin{align}\label{eq:GTAw} 
			\hspace{-0.3mm} w_i^{t+1} &= w_i^t -\gamma\sum_{j\in \mathcal{N}_i}\!\! {\ell}_{ij} (w_j^{t} - \delta d_j^t) -  \gamma s_i^t + \delta(d_i^{t+1} - d_i^t)
			\\[.5em]%
			s_i^{t+1} &= s_i^t - \! \gamma\!\sum_{j\in \mathcal{N}_i} \ell_{ij} s_j^{t} + \frac{2}{\delta}(f_i(w_i^{t+1})d_i^{t+1} -f_i(w_i^t )d_i^t)\! 
			\end{align}
		\end{subequations}
		\EndFor
	\end{algorithmic}
	\caption{\algo/ (agent $i$)}
	\label{table:algorithm}
\end{algorithm}

The local gradient estimates generated by $\frac{f_i(w_i^t)d_i^t}{2\delta}$ are suitably interlaced with (i) the term $\sum_{j\in \mathcal{N}_i} {\ell}_{ij} (w_j^{t} - \delta d_j^t)$ to force consensus among the quantities $w_i^t - \delta d_i^t$ and (ii) a tracking mechanism to reconstruct the (estimated) global gradient.
In detail, the consensus step is performed by using the entries $\ell_{ij}$ the $(i,j)$-entry of the Laplacian matrix $\mathcal{L}$ associated to the graph $\mathcal{G}$, while the tracking mechanism is implemented by equipping each agent $i$ with an auxiliary variable $s_i^t \in \R^{\xD}$.
In this algorithm, agents exchange with their neighbors the information $\col(w_i^t - \delta d_i^t, s_i^t)$ involving $2 \xD$ components.
\begin{remark}
	The main distinctive features of our method are as follows.
	First, errors in our method are given by third-order residuals as
	opposed to second-order ones in finite-difference methods.
	Second, gradient estimation is based on a single-function query per agent. This
	could be advantageous in scenarios in which multiple queries per agent are
	expensive or even not allowed.
	Third and final, the estimation policy updates, and so the convergence
	guarantees, are purely deterministic. \oprocend
\end{remark}

The next theorem formalizes the convergence properties of \algo/.
To this end, for all $i \in \until{N}$ and given $r > 0$, we define the set $\cD_{i,r} \subset \R^{2\xD}$ as 
\begin{align*}
	\cD_{i,r} \! := \! \left\{\col(w_i,s_i) \in \R^{2\xD} \! \mid \! \norm{w_i -
	x^\star} \leq r, s_i = 2f_i(w_i)/\delta\right\}\!.
\end{align*}
\begin{theorem}
	\label{th:MainResult}
	Consider~\eqref{eq:GTA} and let Assumptions~\ref{ass:connectivity},~\ref{ass:strong_convexity}, and~\ref{ass:lipschitz} hold.
	Then, for any $r, \bar{\rho} > 0$, there exist
	$\gamma^\star, \delta^\star, k_1 > 0$, $\epsilon \in (0,\bar{\rho}/2)$, and
	$k_2 \ge (\bar{\rho}/2 - \epsilon)$ such that, for any
	$\gamma \in (0,\gamma^\star)$, $\delta \in (0,\delta^\star)$, and $\col(w_i^0,s_i^0) \in \cD_{i,r}$ for all
	$i \in \until{N}$, the trajectories of~\eqref{eq:GTA} are bounded and satisfy
	\begin{align}\label{eq:result_theorem}
		\norm{w_i^t - x^\star} \leq \bar{\rho},
	\end{align}
	for all $i \in \until{N}$ and $t \ge t^\star :=\! -\frac{1}{\gamma k_1}\ln((\bar{\rho}/2 \! - \!\epsilon)/k_2)$, i.e., the convergence to the set $\{w_i \in \R^\xD \mid \norm{w_i^t \! - \! x^\star} \leq \bar{\rho}\}$ is linear.\oprocend
\end{theorem}

The proof of Theorem \ref{th:MainResult} is provided in Section~\ref{sec:proof_of_theorem}.
Theorem~\ref{th:MainResult} provides a semi-global, practical exponential-stability result restricted to the set $\cD_r := \cD_{1,r} \times \dots \times \cD_{N,r}$. %
Indeed, it is semi-global because the parameters $\gamma^\star$ and $\delta^\star$ depend on the initial radius $r$ and it is practical because they also depend on the arbitrary small final radius $\bar{\rho} > 0$.

\section{\algo/: \\ Stability Analysis}
\label{sec:MainResult}
\begin{figure}[H]
	\centering
	\includegraphics[scale=.8]{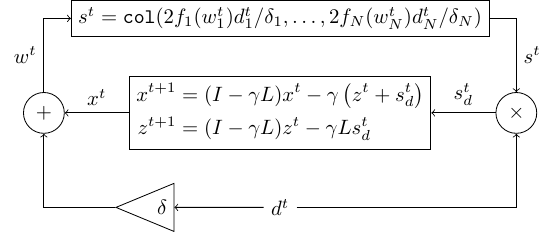}
	\caption{Block scheme of the proposed \algo/ algorithm in the $(x,z)$ coordinates.
      }
	\label{fig:ProposedSolution}
\end{figure}

In this section, we analyze \algo/.
Assumptions~\ref{ass:connectivity},~\ref{ass:strong_convexity}, and~\ref{ass:lipschitz} hold throughout the whole section.
First, let the coordinates $x_i^t, z_i^t \in \R^{\xD}$ be defined as 
\begin{align}\label{eq:x_i_z_i}
	x_i^t := w_i^t -\delta d_i^t, \quad z_i^t := s_i^t - \frac{2f_i(w_i^t)d_i^t}{\delta},
\end{align}
for all $i\in\until{N}$, which allow us to rewrite~\eqref{eq:GTA} as 
\begin{subequations}\label{eq:GTA_x_z}
	\begin{align}
	\hspace{-0.3mm} x_i^{t+1} &= x_i^t -\gamma \bigg(\sum_{j\in \mathcal{N}_i} {\ell}_{ij} x_j^{t} +  z_i^t +\frac{2 f_i(x_i^t+\delta d_i^t)d_i^t}{\delta}\bigg)
	\\%
	z_i^{t+1} &= z_i^t - \gamma\sum_{j\in \mathcal{N}_i} \ell_{ij}\left(z_j^{t}+\frac{2  f_j(x_j^t + \delta d_j^t) d_j^t}{\delta}\right).
	\end{align}
\end{subequations}
\begin{remark}\label{rem:algorithm_CGT}
	The new variables $x_i^t$ and $z_i^t$ allow us to interpret \algo/ as an approximated discrete-time version of the continuous-time gradient tracking method proposed in~\cite{CARNEVALE2023110726}.
	Indeed, if the case gradients were available, the distributed algorithm~\eqref{eq:GTA_x_z} would become%
	\begin{subequations}\label{eq:ClassicGTCausal}
		\begin{align}
			x_i^{t+1} &={x_i^t} -{\gamma}\sum_{j\in \mathcal{N}_i} {\ell}_{ij}x_j^{t} - \gamma \left(z_i^t +\nabla f_i(x_i^t)\right) \\ 
			z_i^{t+1} &={z_i^t} -{\gamma}\sum_{j\in \mathcal{N}_i} {\ell}_{ij}\left(z_j^{t}+\nabla f_j(x_j^t)\right).
		\end{align}
	\end{subequations}
	The analysis required to prove Theorem~\ref{th:MainResult} will also require the investigation of the convergence properties of~\eqref{eq:ClassicGTCausal} (see Lemma~\ref{lemma:Average} and Remark~\ref{rem:per_se_result}). \oprocend
\end{remark}
Then, we aggregate the local updates in~\eqref{eq:GTA_x_z} obtaining the compact algorithm description
\begin{subequations}\label{eq:ESGTcompact}
\begin{align}
		{x}^{t+1} &=\,x^t - \gamma \left(L {x}^t + z^t + f_d(t,x^t)\right)
		\\
		{z}^{t+1} & =\, z^t - \gamma \left(L {z}^t + L f_d(t,x^t)\right),
	\end{align}
\end{subequations}    
where we introduced $L :=\mathcal{L}\otimes I_{\xD}$,
$x^t := \col(x_1^t,\dots,x_N^t)$, $z^t := \mathtt{col}(z_1^t, \dots, z_\nAg^t)$, $d^t := \mathtt{col}(d_1^t, \dots, d_\nAg^t)$, %
and the function $f_d: \N \times \R^{\nAg\xD}\to \R^{\nAg\xD}$ defined as 
\begin{align}
	f_d (t,x) 
	:= 
	\begin{bmatrix}
		2f_1 (x_1 + \delta d_1^t)d_1^t/\delta
		\\
		 \vdots
		 \\
		 2f_\nAg(x_\nAg +\delta d_\nAg^t)d_\nAg^t/\delta
		\end{bmatrix},\label{eq:f_d}
\end{align}
where we decomposed $x$ according to $x := \col(x_1,\dots,x_N)$ with $x_i \in \R^{\xD}$ for all $i \in \until{\nAg}$.
We point out that system~\eqref{eq:ESGTcompact} can be conceived as an extremum seeking scheme
with output map $f(x+\delta d^t)$, see also Fig.~\ref{fig:ProposedSolution}.
We now give an overview of the main steps of the stability analysis carried out to prove
Theorem~\ref{th:MainResult}:
\begin{enumerate} [label=(\roman*)]
\item We perform two change of variables to describe the
  dynamics~\eqref{eq:ESGTcompact} in terms of the mean value (over the agents) $\tilde{z}_{\text{avg}}$ of
  $z$ and the orthogonal part $\tilde{z}_\perp$ (both in error coordinates). %
  Then, by relying on averaging theory (see~\cite{bai1988averaging} for discrete-time systems or~\cite[Ch.10]{khalil2002nonlinear} and~\cite{sanders2007averaging} for continuous-time ones), we introduce a suitable auxiliary system named \emph{averaged system}.
	The latter is obtained by averaging the original algorithm dynamics over a
  common period.
 The averaged system is shown to be driven by the cost gradients with
  additive estimation errors.

\item When neglecting these errors, the averaged system corresponds to an
    equivalent form of \eqref{eq:ClassicGTCausal}.
  Based on this observation, we rely on existing stability properties of the
    continuous gradient tracking to demonstrate that the trajectories of the
    averaged system exponentially converge to an arbitrarily small neighborhood of
    $\col(\1 x^\star, \tilde{z}^{\text{eq}}_\perp)$ for some $\tilde{z}^{\text{eq}}_\perp$
    arising from the analysis.
\item Finally, we prove Theorem~\ref{th:MainResult} by exploiting the steps
    above and by using averaging theory to show the closeness between the trajectories
    of~\eqref{eq:ESGTcompact} and those of its averaged system.
\end{enumerate}

Step (i) is performed in Section~\ref{sec:ChangeOfCoordinate}, step (ii)
is carried out in
Section~\ref{sec:averaging}, while Section~\ref{sec:proof_of_theorem}
is devoted to step (iii).

\subsection{Coordinate changes and averaged system}
\label{sec:ChangeOfCoordinate}

We start by introducing a change of coordinates to highlight the error dynamics of~\eqref{eq:ESGTcompact} with respect to $\col(x^\star,-G(\1x^\star))$, i.e., the point in which each $x_i^t$ coincides with the optimal problem solution $x^\star$ and perfect tracking is achieved via $z_i^t$ (see~\cite{CARNEVALE2023110726}).
To this end, let $G: \R^{\nAg\xD} \to \R^{\nAg\xD}$ be defined as 
\begin{align}
 G(x) :=&\, \mathtt{col}(\nabla f_1(x_1),\dots,\nabla f_\nAg (x_\nAg)).\label{eq:G}
\end{align}
Then, let the error coordinates $\tilde{x}, \tilde{z} \in \R^{\nAg\xD}$ be defined as
\begin{align}
	\tilde{x} := x - \1 x^\star, \quad \tilde{z} = z + G(\1 x^\star),\label{eq:change_error}
\end{align} %
and let us introduce $\phi_{xz}: \N \times \R^{\nAg\xD} \times \R^{\nAg\xD} \to \R^{2\nAg\xD}$ as
\begin{align*}
\phi_{xz}(t,\tilde{x},\tilde{z}) \! \!:=\!\! \begin{bmatrix}
- L \tilde{x}  - \tilde{z}  -  f_d(t,\tilde{x} +\1 x^\star) + G(\1 x^\star)
\\
- L \tilde{z}  - L (f_d(t,\tilde{x} +\1 x^\star) - G(\1 x^\star))
\end{bmatrix}\!.
\end{align*}
Then, by using the new coordinates, we rewrite~\eqref{eq:ESGTcompact} as
\begin{equation}
	\label{eq:xz_tilde}
	\begin{bmatrix}\tilde{x}^{t+1}\\
		\tilde{z}^{t+1}
	\end{bmatrix}
	= 	\begin{bmatrix}\tilde{x}^t\\
		\tilde{z}^{t}
	\end{bmatrix} +\gamma \phi_{xz}(t,\tilde{x}^t,\tilde{z}^t),%
\end{equation}
where we have used the property $L\1 = 0$. %
As in~\cite{CARNEVALE2023110726}, we take advantage of the initialization $z_i^0 = 0$ for all $i \in \until{N}$.
To this end, we introduce the novel coordinates $\tilde{z}_{\text{m}} \in \R^{\xD}$ and $\tilde{z}_\perp \in \R^{(\nAg-1)\xD}$ representing the average of the variables $\tilde{z}_i^t$ (over the agents) and the orthogonal counterpart, namely
\begin{align}\label{eq:change}
  \begin{bmatrix}
    \tzm
    \\
    \tilde{z}_\perp
  \end{bmatrix} := \begin{bmatrix}
    \frac{\1^\top}{N}
    \\
    R^\top
  \end{bmatrix}\tilde{z},
\end{align}
where we introduced the matrix $R \in \R^{\nAg\xD \times (\nAg - 1)\xD}$ such that $R^\top \1 = 0$, $R^\top R = I$.
Further, given $n_{\xi} := (2\nAg-1)\xD$, let us also introduce $\xi \in \R^{n_\xi}$ defined as
\begin{align}
	\xi := \begin{bmatrix}
    \tilde{x}^\top
    &
    \tilde{z}_\perp^\top
  \end{bmatrix}^\top.
  \label{eq:change_xi}
\end{align}
Since $\1^\top L = 0$ in light of Assumption~\ref{ass:connectivity} and $\1^\top G(\1 x^\star) = \sum_{i=1}^N \nabla f_i(x^\star) =0$, we use~\eqref{eq:change} and~\eqref{eq:change_xi} to rewrite~\eqref{eq:xz_tilde} as 
\begin{subequations}\label{eq:xi_barz}
  \begin{align}
    \xi^{t+1} &= \xi^t + \gamma \phi_{\xi}(t,\col(\tilde{x}^t,\1\tzm^t + R\tilde{z}_\perp^t))\label{eq:xi_dynamics_with_barz}
    \\
    \tzm^{t+1} &= \tzm^t,\label{eq:barz_dynamics}
  \end{align}
\end{subequations}
where we introduced $\phi_\xi: \N \times \R^{2\nAg\xD} \to \R^{n_{\xi}}$ defined as
\begin{align}
  \phi_\xi(t,\col(\tilde{x},\tilde{z})) := \begin{bmatrix} I& 0\\ 
    0& R\end{bmatrix}\phi_{xz}(t,\tilde{x},\tilde{z}).\label{eq:phi_xi}
\end{align}
The equation~\eqref{eq:barz_dynamics} implies that $\tzm^t = \tzm^0$ for all $t \in \N$.
Then, since $z_i^0 = 0$ for all $i \in \until{N}$ and
$\1^\top G(\1x^\star) =0$, it holds $\tzm^0 = 0$ which allows us to ignore~\eqref{eq:barz_dynamics} and rewrite~\eqref{eq:xi_barz} according to the equivalent, reduced system 
\begin{align}\label{eq:xi}
  \xi^{t+1} = \xi^t + \gamma \phi(t,\xi^t),
\end{align}
where $\phi: \N \times \R^{n_\xi} \to \R^{n_\xi}$ is introduced to compactly describe system~\eqref{eq:xi_dynamics_with_barz} with $\tzm^t = 0$ for all $t \in \N$, i.e., $\phi$ is defined as
\begin{align}
\phi(t,\xi) &:= \phi_{\xi}(t,\col(\tilde{x},R\tilde{z}_\perp))%
\label{eq:phi}\\
&\stackrel{(a)}{=} \!\begin{bmatrix}
	I& 0
	\\
	0& R
\end{bmatrix}\!\!\begin{bmatrix}
	- L \tilde{x}  - R\tzperp  -  f_d(t,\tilde{x} +\1 x^\star) + G(\1 x^\star)
	\\
	- LR\tzperp  - L (f_d(t,\tilde{x} +\1 x^\star) - G(\1 x^\star))
	\end{bmatrix}\!\!,\notag
\end{align}
where in $(a)$ we used the definition of $\phi_\xi$ (cf.~\eqref{eq:phi_xi}).

	We resort to the averaging theory \cite{bai1988averaging} to analyze the time-varying system~\eqref{eq:xi}. 
	As customary when this tool is employed, we introduce an auxiliary scheme typically named \emph{averaged system}, obtained by averaging the time-varying vector field $\phi(t,\xi)$ over $\per$ samples by \emph{freezing} the system state. In detail, the averaged system associated with \eqref{eq:xi} is
\begin{equation}\label{eq:xia_implicit}
	 \xi\av^{t+1} =  \xi\av^t
	+\gamma\phi\av(\xi\av^t), \quad \quad \quad \text{(averaged system)} 
\end{equation} 
with $\xi\av^0 = \xi_0$, where $\phi\av: \R^{n_\xi} \to \R^{n_\xi}$ is defined as
\begin{align}
	\phi\av(\xi) := \frac{1}{\per}\sum_{k=t + 1}^{t + \per} \phi(k,\xi), \quad \text{for any } t\ge0.\label{eq:phi_av}
\end{align}
Notice that, since $\phi(k,\xi)$ is periodic with period $\per$ on the first argument, the function $\phi\av(\xi)$ (cf.~\eqref{eq:phi_av}) is time-invariant.

The next lemma shows that the gradient of each $f_i$ can be approximated by averaging each block of $f_d$ (cf.~\eqref{eq:f_d}) over the period $\per$ on the first argument.
For this reason, the next lemma will be useful to explicitly write $\phi\av$ (cf.~\eqref{eq:phi_av}) in terms of $G(x)$, (i.e., the stack of the gradients $\nabla f_i(x_i)$, see~\eqref{eq:G}).
\begin{lemma}[Gradient estimation]
\label{lemma:ESApprox}
For all $i \in \until{N}$, there exists $\ell_i:\mathbb{R}^{\xD}\to\mathbb{R}^{\xD}$ such that, for any given $x_i \in \R^{\xD}$ and all $t \in \N$, it holds
\begin{align}
	\dfrac{2}{\delta\per}\sum_{k=t+1}^{t + \per} f_i (x_i+\delta d_i^k)d_i^k = \nabla  f_i(x_i) + \delta^2	 \ell_i(x_i).\label{eq:result_gradient_estimation}
\end{align}
Moreover, given any compact set $\mathcal{S}_i \subset \R^\xD$, if $\delta \in (0,1]$, there exists $L_{i,\mathcal{S}_i} > 0$ such that 
\begin{align}
	\norm{\ell_i(x_i)} \leq L_{i,\mathcal{S}_i},\label{eq:ell_bound}
\end{align}
for all $x_i \in \mathcal{S}_i$ and $i \in \until{N}$. \oprocend

\end{lemma}
The proof of Lemma \ref{lemma:ESApprox} is in Appendix~\ref{app:ProofLemmaESApprox}.

	We note that the quantity on the right-hand side of~\eqref{eq:result_gradient_estimation} is time-invariant because  $d_i^k$ is periodic, with period $\per$, and the left-hand side of~\eqref{eq:result_gradient_estimation} is averaged over the period $\per$.
Now, let us introduce the function $\ell: \R^{\nAg\xD} \to \R^{\nAg\xD}$ stacking all the approximation errors $\ell_i(x_i)$ used in~\eqref{eq:result_gradient_estimation}, namely
\begin{align}
	\ell(x) := \mathtt{col}\left( \ell_1(x_1), \dots,  \ell_\nAg(x_\nAg)\right).\label{eq:ell}
\end{align}
Then, by using~\eqref{eq:result_gradient_estimation}, the definitions of $f_d$ in~\eqref{eq:f_d}, $G(x)$ as the stack of the gradients $\nabla f_i(x_i)$ in~\eqref{eq:G}, and $\ell(x)$ as the stack of the approximation errors $\ell_i(x_i)$ in~\eqref{eq:ell}, it holds  
\begin{align}
	\frac{1}{\per}\sum_{k=t+1}^{t+\per} f_d(k,x) = G(x) + \ell(x),\label{eq:f_d_G}
\end{align}
for all $x \in \R^{\nAg\xD}$ and $t \in \N$.
Hence, by introducing $\xi\av := \col(\tilde{x}\av,\tilde{z}_{\perp,\av})$ and by combining the definition of $\phi\av$ in~\eqref{eq:phi_av}, $\phi$ in~\eqref{eq:phi}, and~\eqref{eq:f_d_G}, we obtain
\begin{align}
	\phi\av(\xi\av) \! &= \! \begin{bmatrix}
		- L \tilde{x}\av  -  R\tilde{z}_{\perp,\av}  -  G(\tilde{x}\av +\1 x^\star) + G(\1 x^\star)
		\\
		- R^\top L R\tilde{z}_{\perp,\av}  \! - \!  R^\top L (G(\tilde{x}\av +\1 x^\star) \! - \! G(\1 x^\star))
		\end{bmatrix}
		\notag\\
		&\hspace{.4cm}
		+ \begin{bmatrix}
			- \ell(\tilde{x}\av + \1x^\star)
			\\
			-R^\top L\ell(\tilde{x}\av + \1x^\star)
		\end{bmatrix}.\label{eq:phi_av_explicit}
\end{align}
We define 
\begin{align*}
	\phi_{GT}(\xi\av) \! &:= \! 
	\begin{bmatrix}
		- L \tilde{x}\av  -  R\tilde{z}_{\perp,\av}  -  G(\tilde{x}\av +\1 x^\star) + G(\1 x^\star)
		\\
		\! - \! R^\top L R\tilde{z}_{\perp,\av}  \! - \!  R^\top L (G(\tilde{x}\av \! + \! \1 x^\star) \! - \! G(\1 x^\star))
	\end{bmatrix}
	\\
	B \! &:= \!
	\begin{bmatrix}
		-I\\
		-R^\top L
	\end{bmatrix}
,\quad 
u(\xi\av):= \ell(\tilde{x}\av + \1 x^\star)
\end{align*}
to short~\eqref{eq:phi_av_explicit} as
\begin{align*}
	\phi\av(\xi\av) = \phi_{GT}(\xi\av) + \delta^2 B u(\xi\av),
\end{align*}
which, in turn, allows us to rewrite~\eqref{eq:xia_implicit} as 
\begin{equation}\label{eq:xia}
	\xi\av^{t+1} =  \xi\av^t + \gamma\phi_{GT}(\xi\av^t) + \gamma \delta^2 B u(\xi\av^t) \quad \text{(averaged system)}.
\end{equation} 

\subsection{Averaged System Analysis}
\label{sec:averaging}

In this subsection, we analyze the averaged system~\eqref{eq:xia}. %
To this end, we first consider an additional nominal system in which the term $\gamma \delta^2 B u(\xi\av^t)$ (i.e., the term describing the gradients' estimation error) is neglected.
Then, by using such a nominal system analysis as a building block, we provide the result concerning system~\eqref{eq:xia}.
Therefore, we start by studying
\begin{equation}\label{eq:GTAUnperturbed}
\xi\av^{t+1} = \xi\av^t + \gamma \phi_{GT}(\xi\av^t),
\end{equation}
which corresponds to system~\eqref{eq:xia} in the case of $\gamma \delta^2 B u(\xi\av^t) = 0$.
The next lemma proves the global exponential stability of the origin
for~\eqref{eq:GTAUnperturbed}. %
\begin{lemma}
	\label{lemma:AVunperturbed} %
	There exist $P = P^\top \in \R^{n_\xi \times n_\xi }$ and $a_1, a_2, c_1, \gamma_0  > 0$ such that, for any $\gamma \in (0,\gamma_0)$, along the trajectories of~\eqref{eq:GTAUnperturbed} it holds 
	\begin{subequations}\label{eq:V_AVunperturbed}
		\begin{align}
			a_1I \leq P &\leq a_2 I\label{eq:V_bounds}
			\\
			{\xi\av^{t+1}}^\top P \xi\av^{t+1} - {\xi\av^t}^\top P\xi\av^t  &\leq - \gamma c_1\norm{\xi\av^t}^2,\label{eq:V_negative}
		\end{align}
	\end{subequations}
	for all $\xi\av^t \in \R^{n_\xi}$.\oprocend
\end{lemma}
The proof of  Lemma~\ref{lemma:AVunperturbed} is in Appendix~\ref{sec:ProofLemmaAVunperturbed}.

\begin{remark}\label{rem:per_se_result}
  Notice that Lemma~\ref{lemma:AVunperturbed} proves that algorithm~\eqref{eq:ClassicGTCausal} linearly
  converges to the minimizer of \eqref{eq:TotalCost}, since
  \eqref{eq:GTAUnperturbed} is an equivalent formulation of~\eqref{eq:ClassicGTCausal}.\oprocend
\end{remark}

With this result at hand, we analyze the impact of $u(\cdot)$ thus obtaining the stability properties of the averaged system~\eqref{eq:xia}.

\begin{lemma}\label{lemma:Average}
  Consider the averaged system~\eqref{eq:xia}.
  Then, for any $r_{\xi} > 0$ and
  $\rho \in (0,r_{\xi})$, there exist $c_3 \in (0,c_1)$ and $\delta^\star_1 \in (0,1]$
  such that, for any $\gamma \in \min\{\gamma_0,1\}$, $\delta \in (0,\delta^\star_1)$, and $\norm{\xi\av^0} \leq r_{\xi}$, it holds
  \\
(i) $\xi\av^t \in \mathcal{B}_{\sqrt{a_2/a_1} r_\xi}$ for all $t \in \N$,
\\
(ii)
  \begin{align}
\|\xi\av^t\| \le \sqrt{a_2/a_1} \exp\left(-t \gamma c_3 \right)  \|\xi\av^0\|, \label{eq:exponential_lemma}
  \end{align}
for all $\norm{\xi\av^t} \ge  \rho$.\oprocend
\end{lemma}
The proof of Lemma~\ref{lemma:Average} is in Appendix~\ref{sec:ProofLemmaAverage}.

\begin{remark}
	The result of Lemma~\ref{lemma:Average} only involves the averaged system~\eqref{eq:xia}.
	In the next section, such a result will be used as a building block to study the original dynamics, i.e., system~\eqref{eq:xi} and, thus, to conclude the proof of Theorem~\ref{th:MainResult}.
	However, it is a per se result that can be used to show robust stability for the distributed algorithm~\eqref{eq:ClassicGTCausal}.\oprocend
\end{remark}

\subsection{Proof of Theorem \ref{th:MainResult}}
\label{sec:proof_of_theorem}

Despite averaging tools for discrete-time systems are already present in the literature, see \cite{bai1988averaging} for example, we got the inspiration from continuous-time averaging \cite[Ch.~10]{khalil2002nonlinear} and~\cite{sanders2007averaging} for elaborating the proof of Theorem \ref{th:MainResult} to make clear how $\gamma$ affects the closeness of the trajectories of \eqref{eq:xi} and \eqref{eq:xia}.
Since Assumptions~\ref{ass:connectivity},~\ref{ass:strong_convexity}, and~\ref{ass:lipschitz} hold, we apply Lemma~\ref{lemma:AVunperturbed} to claim that there exist $P = P^\top \in \R^{n_\xi \times n_\xi }$ and $a_1, a_2, c_1,\gamma_0 > 0$ such that, if $\gamma \in (0,\gamma_0)$, the conditions~\eqref{eq:V_AVunperturbed} are satisfied.
Then, we evaluate the norm of the initial conditions of system~\eqref{eq:xi} and~\eqref{eq:xia}, i.e., $\norm{\xi^0} = \norm{\xi^0\av}$.
By using the definition of $\xi$ (cf.~\eqref{eq:change_xi}), the changes of variables~\eqref{eq:change_error} and~\eqref{eq:change}, and the triangle inequality, we get 
\begin{align*}
	\norm{\xi^0} &\leq \norm{x^0 - \1x^\star} + \norm{R^\top(z^0 + G(\1x^\star))} 
	\\
	&\hspace{.4cm}
	+ \norm{\frac{\1^\top}{N}(z^0 + G(\1 x^\star))}
	\stackrel{(a)}{\leq} r\sqrt{\nAg} + \norm{R}\norm{G(\1 x^\star)},
\end{align*} 
where in $(a)$ we combine the initialization $\norm{x_i^0 - x^\star} \leq r$ and $z_i^0 = 0$ for all $i \in \until{N}$ with the fact that $\1^\top G(\1 x^\star) = \sum_{i=1}^N f_i(x^\star) = 0$.
Hence, by defining $r_{\xi} := r\sqrt{\nAg} + \norm{R}\norm{G(\1 x^\star)}$, we claim that 
\begin{align*}
	\norm{\xi^0} = \norm{\xi\av^0} \leq r_\xi.
\end{align*}
Once the initial distance from the origin has been evaluated, we choose any $\bar{\rho} > 0$, set $c_2 := \sqrt{a_2/a_1}$, and choose any $\epsilon \in (0,\bar{\rho}(1+c_2)/2)$.
Then, we pick $\rho \in (0,(\bar{\rho}/2 - (1 + c_2)\epsilon)/c_2)$, $c_1^\prime \in (0,c_1)$, and use the matrix $P$ satisfying~\eqref{eq:V_AVunperturbed} to apply Lemma~\ref{lemma:Average}. 
Specifically, we claim that there exist $c_3 >0$, $\delta^\star_1 \in (0,1]$, and $\gamma_0 > 0$ such that, for any $\delta \in (0,\delta^\star_1)$ and $\gamma \in (0,\min\{\gamma_0,1\})$, it holds $\xi\av^t \in \mathcal{B}_{c_2 r_\xi}$ for all $t \in \N$ and %
\begin{align}
	\norm{\xi\av^t} \leq c_2\exp(-t\gamma c_3)\norm{\xi\av^0},\label{eq:exponential_proof}
\end{align}
for all $\xi\av^t$ such that $\norm{\xi\av^t} \ge \rho$.
Now, we proceed by finding a bound $\gamma_1 > 0$ such that, for any $\gamma \in (0,\gamma_1)$, we guarantee the $\epsilon$-closeness between the state $\xi\av^t$ of the averaged system~\eqref{eq:xia_implicit} and the one of the original system~\eqref{eq:xi}, i.e., that $\|\xi^t-\xi\av^t\| \leq \epsilon$ holds true for all $t \in \N$.
To this end, 
let us introduce
\begin{align*}
	\upsilon(t,\xi\av) := \sum_{k=0}^{t-1}\left(\phi(k,\xi\av)-\phi\av(\xi\av)\right).
\end{align*}
By using this definition and the one of $\phi\av$ (cf.~\eqref{eq:phi_av}), it holds
\begin{align}
&\upsilon(t+1,\xi\av^{t+1})- \upsilon(t,\xi\av^t) 
	\notag\\
			&= \phi(t,\xi\av^{t+1})-\phi\av(\xi\av^{t+1})+ \upsilon(t,\xi\av^{t+1}) - \upsilon(t,\xi\av^t).\label{eq:upsilon}
\end{align}
Then, let $r_\xi^\prime := c_2r_\xi$ and define $\Delta := \delta\sqrt{\nAg\xD}$. %
Under the assumption of $\xi^t \in %
		\mathcal{B}_{r_{\xi}^\prime + \epsilon}$ for all $t \in \N$ (later verified by a proper selection of $\gamma$), we claim that the arguments of the functions $f_i$ and their derivatives (embedded into the definitions of $\phi(t,\cdot)$ and $\phi\av(\cdot)$ and their derivatives) lie into the compact set $\mathcal{B}_{r^{\prime}_\xi + \epsilon + \Delta}$.
		Thus, since the functions $f_i$ and its derivatives are continuous (cf. Assumption~\ref{ass:lipschitz}) and the functions $\phi(\cdot,\cdot)$ and $\nu(\cdot,\cdot)$ are periodic in the first argument, we define 
		\begin{align}
			L_\phi := \sup_{\substack{
			\xi \in \mathcal{B}_{r_{\xi}^\prime + \epsilon}
			\notag\\
			t \in [0,\per]}}\bigg\{&\norm{\phi(t,\xi)}, \norm{\phi\av(\xi)}, \norm{\frac{\partial \phi(t,\xi)}{\partial \xi}}, 
			\\
			&\norm{\frac{\partial \phi\av(\xi)}{\partial \xi}}, \norm{\frac{\partial \nu(t,\xi)}{\partial \xi}}\bigg\}.\label{eq:L_phi}
		\end{align}
		Consequently, for all $\xi, \xi^\prime \in \mathcal{B}_{r_\xi^\prime + \epsilon}$ and $t \in \N$, it holds
		\begin{subequations}\label{eq:bounds_phi}
		\begin{align}
		\|\upsilon(t,\xi)\| &\leq 2 L_\phi \per\label{eq:bound_upsilon}
		\\
		\|\phi(t,\xi)-\phi(t,\xi^\prime)\| &\leq  L_\phi \|\xi-\xi^\prime\|
		\\
		\|\phi\av(\xi)-\phi\av(\xi^\prime)\| &\leq  L_\phi \|\xi-\xi^\prime\|
		\\
		\| \upsilon(t,\xi) - \upsilon(k,\xi^\prime)\| &\leq 2L_\phi \per \|\xi-\xi^\prime\|
		\\
		\norm{\phi\av(\xi)} &\leq L_\phi.
		\end{align}
		\end{subequations}
		Let us introduce $\zeta^t \in \R^{n_\xi}$ defined as
		\begin{align}
			\zeta^t := \xi\av^t+\gamma \upsilon(t,\xi\av^t).\label{eq:zeta_definition}
		\end{align} %
		Then, it holds
		\begin{align*}
			\xi^t-\zeta^t &=  \sum_{k=0}^{t-1} (\xi^{k+1}-\xi^k)-(\zeta^{k+1}-\zeta^k),
		\end{align*}
	add $\pm\gamma \sum_{k=0}^{t-1}(\phi(k,\zeta^k)  + \phi(k,\xi\av^k)),
	$
	and use~\eqref{eq:upsilon} to get
	\begin{align*}
		&\xi^t-\zeta^t =  \gamma \sum_{k=0}^{t-1}(\phi(k,\xi^k)-\phi(k,\zeta^k))
		\\
		&+\gamma\sum_{k=0}^{t-1}(\phi(k,\zeta^k)-\phi(k,\xi\av^k))
		\\
		&\hspace{.4cm}
		-\gamma\sum_{k=0}^{t-1}(\phi(k,\xi\av^{k+1})-\phi(k,\xi\av^k))
		\\
		&+\gamma\sum_{k=0}^{t-1}(\phi\av(\xi\av^{k+1})-\phi\av(\xi\av^k))
		\\
		&\hspace{.4cm}
		-\gamma\sum_{k=0}^{t-1} ( \upsilon(k,\xi\av^{k+1}) - \upsilon(k,\xi\av^k)).
	\end{align*}
	Use \eqref{eq:xi}, \eqref{eq:xia}, and~\eqref{eq:bounds_phi} to bound 
	\begin{align}
		\|\xi^t-\zeta^t\| 
		&\le  \gamma L_\phi \sum_{k=0}^{t-1}\|\xi^k-\zeta^k\|+\gamma^2 L_\phi^2 2 \left( 1 + 2  \per \right) t.\label{eq:xi_minus_zeta}
	\end{align}
	Apply the discrete Gronwall inequality (see~\cite{Popenda1983Discrete,holte2009discrete})
	and
	\[
	\sum_{k=0}^{t-1}  \gamma L_\phi k  \exp\left(- \gamma L_\phi k \right) \le \sum_{k=0}^{\infty}  \gamma L_\phi k  \exp\left(- \gamma L_\phi k \right) = 1
	\]
	to further bound~\eqref{eq:xi_minus_zeta} as
	\begin{align}
		\|\xi^t-\zeta^t\| 
		&\le \gamma^2 L_\phi^2 2 \left( 1 + 2  \per \right) t
		\notag\\
		&\hspace{.4cm} +  
		\gamma L_\phi 2 \left( 1 + 2  \per \right)  \exp\left(\gamma L_\phi t \right).\label{eq:xi_minus_zeta_2}
	\end{align}
	The definition of $\zeta^t$~\eqref{eq:zeta_definition} and the triangle inequality lead to
	\begin{align}
			\|\xi^t-\xi\av^t\| &\leq \norm{\xi^t - \zeta^t} + \gamma\norm{v(t,\xi\av^t)}
			\notag\\
			&\stackrel{(a)}{\leq}
			\gamma^2 L_\phi^2 2 \left( 1 + 2  \per \right) t   +\gamma 2L_\phi \per
			\notag\\
			&\hspace{.4cm}
		+\gamma L_\phi 2 \left( 1 + 2  \per \right)  \exp\left(\gamma L_\phi t \right).\label{eq:xi_minus_xi_avg}
	\end{align}
	where $(a)$ uses~\eqref{eq:xi_minus_zeta_2} to bound the first term and~\eqref{eq:bound_upsilon} to bound the second one.
	Then, set $\theta^\star \in \real$ such that
	\begin{align}
		\theta^\star \ge -\tfrac{1}{c_3}\ln\left(\tfrac{(\bar{\rho}/2 - \epsilon)/c_2}{c_2 r_\xi}\right).\label{eq:theta_star}
	\end{align}
	Let
		$\gamma_2 := \tfrac{\epsilon/(3c_2)}{L_\phi^2 2 \left( 1 + 2  \per \right) \theta^\star}$, %
		$\gamma_3 := \tfrac{\epsilon/(3c_2)}{2L_\phi (1+2\per)\exp(L_\phi \theta^\star)}
		$,%
		$\gamma_4 := \tfrac{\epsilon/(3c_2)}{2 L_\phi \per}$,
	$\gamma_1 := \min\{\gamma_0,\gamma_2, \gamma_3,\gamma_4,1\}$.
	Pick $\gamma \in (0,\gamma_1)$ such that
	$t^\star := \frac{\theta^\star}{\gamma}\in\natural$. 
	This can be done without loss of generality since $\theta^\star$ is a design parameter.
	Then,  we bound~\eqref{eq:xi_minus_xi_avg} as %
	\begin{align}
		\|\xi^t-\xi\av^t\| %
		\leq \frac{\epsilon}{c_2},
		\label{eq:closeness}
	\end{align}
	for all $t \in \{0,\dots,t^\star\}$.
	As a consequence, since $\xi\av^t \in \mathcal{B}_{r_\xi^\prime}$ for all $t \in \N$, it holds ${\xi}^t \in %
	\mathcal{B}_{r_\xi^\prime + \epsilon}$ for all $t \in \{0,\dots,t^\star\}$, i.e., we have verified that the bounds~\eqref{eq:bounds_phi} can be used into the interval $\{0,\dots,t^\star\}$.
	Moreover, the exponential law~\eqref{eq:exponential_proof} and the expression of $\theta^\star$ (cf.~\eqref{eq:theta_star}) ensure that it holds
	\begin{align}
			\|\xi\av^{t}\| %
			\leq (\bar{\rho}/2 - \epsilon)/c_2,\label{eq:exponential_tstar}
	\end{align}
		for all $t \ge t^\star$.
		Now, by using the triangle inequality, we write
		\begin{align}
			\|\xi^{t^\star}\| \le \|\xi^{t^\star}-\xi\av^{t^\star}\| + \|\xi\av^{t^\star}\| 
			\stackrel{(a)}{\le} \frac{\bar{\rho}}{2c_2},\label{eq:ball}
		\end{align}
		where in $(a)$ we combined~\eqref{eq:closeness} and~\eqref{eq:exponential_tstar}.
	    The inequality~\eqref{eq:ball} guarantees that $\xi^{t^\star} \in \mathcal{B}_{\frac{\bar{\rho}}{2c_2}}$, hence we proved that the trajectories of~\eqref{eq:xia} enters into $\mathcal{B}_{\frac{\bar{\rho}}{2c_2}}$ with linear rate.
		Next, in order to show that $\xi^t \in \mathcal{B}_{\bar{\rho}/2}$ for all $t \ge t^\star$, we divide the set of natural numbers in intervals as 
			$\mathbb{N} = \{0,\dots,t^\star\} \cup \{t^\star,\dots,2t^\star\} \cup \dots$.
		Define $\psi\av(k+ t^\star,\xi^{t^\star})$ as the solution to~\eqref{eq:xia} for $\xi\av^0 = \xi^{t^\star}$ and $k \in  \{0,\dots,t^\star\}$.
		Thus, at the beginning of each interval $\{t^\star,\dots,2t^\star\}$, the initial condition of~\eqref{eq:xia} coincides with the one of $\psi\av(k+t^\star,\xi^{t^\star})$ and lies into $\mathcal{B}_{\bar{\rho}/2} \subseteq \mathcal{B}_{r_\xi}$.
		Thus, we apply the same arguments above to guarantee that
		\begin{align*}
			\|\xi^{k+t^\star}-\psi\av(k+t^\star,\xi^{t^\star})\| &\leq \epsilon
			\\
			\psi\av(2t^\star,\xi^{t^\star}) &\in \mathcal{B}_{(\bar{\rho}/2 - \epsilon)/c_2},
		\end{align*}
		for all $\gamma \in (0,\gamma^\star)$ and $k \in  \{0,\dots,t^\star\}$.
		By using Lemma~\ref{lemma:Average}, we guarantee that system~\eqref{eq:xia} cannot escape from $\mathcal{B}_{\bar{\rho}/2 - \epsilon}$, namely $\xi\av^t \in \mathcal{B}_{\bar{\rho}/2 - \epsilon}$ for all $t \ge t^\star$.
		Thus, we get $\xi^t \in \mathcal{B}_{\bar{\rho}/2}$ for all $t \in  \{t^\star,\dots,2t^\star\}$.
		By recursively applying the same arguments above for each interval $\{jt^\star,\dots,(j+1)t^\star\}$ with $j = 2, 3, \dots$ and using $\norm{x_i^t - x^\star} \leq \norm{\xi^t}$ for all $i \in \until{N}$ and $t \in \N$, we get 
		\begin{align}
			\norm{x_i^t - x^\star} \leq \bar{\rho}/2,\label{eq:last_inequality}
		\end{align}
		for all $i \in \until{N}$ and $t \in \N$.
		The change of coordinates~\eqref{eq:x_i_z_i},~\eqref{eq:last_inequality}, and the triangle inequality lead to
		\begin{align}
			\label{eq:boundwi}
			\norm{w_i^t - x^\star} &\leq \frac{\bar{\rho}}{2} + \delta\norm{d_i^t}
			\stackrel{(a)}{\leq} \frac{\bar{\rho}}{2} + \delta\sqrt{\xD},
		\end{align} 
		where in $(a)$ we use the boundedness of the dither signals.			
		The proof follows from~\eqref{eq:boundwi} by setting $\delta^\star := \min\left\{\delta^\star_1,\frac{\bar{\rho}}{2\sqrt{\xD}}\right\}$.

\section{Numerical Computations on Distributed personalized optimization}
\label{sec:Simulation}
To corroborate the theoretical analysis, in this section, we provide numerical
computations for the proposed distributed algorithm on a personalized
optimization framework. %

In several engineering applications, a problem of interest consists of optimizing a performance metric while keeping into account user discomfort terms~\cite{Ye2016Distributed,ospina2021time}.
In these scenarios, the user discomfort term is usually not known in advance but can be only accessed by measurements.
Specifically,
 we associate to each agent $i\in\until{\nAg}$ a cost function in the form
$f_i(w) = w^\top Q_i w + r_i^\top w + \log(\sum_{\ell=1}^\xD a_{i\ell}e^{b_{i\ell} w_\ell})$ 
 with $Q_i = Q_i^\top \in\real^{\xD\times \xD}$, $r_i\in \real^\xD$ and $a_{i\ell},b_{i\ell} > 0$, for all $\ell\in\until{\xD}$.
In the following, we provide different sets of simulations to study in detail the features of \algo/.
In particular, each set consists of Monte Carlo simulations over $20$ randomly generated scenarios in which each one differs from the other in terms of cost functions, communication graphs, and algorithmic variables' initialization.
In particular, unless differently stated, in each trial the agents communicate according
to Erd\H{o}s-R\'{e}nyi random graphs (see, e.g.,~\cite{erdHos1960evolution}) with edge probabilities equal to
$0.2$.
As for the problem parameters, for each trial and all $i \in \until{\nAg}$, we generate each matrix $Q_i$ by pre- and post-multiplying a diagonal matrix (whose diagonal elements are randomly extracted from the interval $[10^{-3}, 5\cdot10^{-3} ]$ with uniform probability) with an orthonormal matrix (randomly generated by extracting its elements from the interval $[0,1]$ with uniform probability) and its transpose, respectively.
Further, for all $i \in \until{N}$ and $\ell \in \until{\xD}$,
we randomly extract the components of $r_i$ within the interval $[- 10^{-2}, 3 \cdot 10^{-2}]$ and the parameters $a_{i\ell},b_{i\ell}$ within the interval $[0,10^{-3}]$ with a uniform probability.
In each set of simulations, we choose the parameters $\tau_{i_p}$ and $\phi_{i_p}$ 
according to the following procedure.
For all $i \in \until{N}$, we take %
	$\phi_{i_p} = \dfrac{\pi}{4}\left(1+(-1)^p\right)$ for all $p = 1,\dots, \xD$,
while $\tau_{i_p}$ have been chosen as the first $\lfloor (\xD+1)/2\rfloor$ elements of odd numbers greater than $3$ since it is possible to show that such a set of frequencies satisfy~\eqref{eq:frequencies}.
Roughly speaking, $\phi_{i,p} = 0$ and $\phi_{i,p} = \pi/2$ are used to create orthogonal functions with the same frequency.
It is important to note that this selection is not unique.
As for the algorithm parameters $\gamma$ and $\delta$, due to the complexity of the dependencies of the bounds $\gamma^\star$ and $\delta^\star$ provided in Theorem~\ref{th:MainResult}, we choose them via a trial-and-error procedure.
In each set of simulations, the performance is evaluated by providing graphical results involving 
  the
  relative errors
  ${|\sum_i f_i(\bar{x}^t)-f(x^\star)|}/{|f(x^\star)|}$
and
  ${\Vert \bar{x}^t- x^\star\Vert}/{\Vert x^\star\Vert}$ achieved along the trials of the Monte Carlo simulations, where $\bar{x}^t:=\frac{1}{N}\sum_{i=1}^N x_i^t$.
  Simulations are performed using DISROPT~\cite{farina2020disropt}, a Python package based on MPI to encode and simulate distributed optimization algorithms.

\subsection{Monte Carlo simulations for large-scale problems}
  
\begin{figure}%
	\centering
	\includegraphics[width=.4975\columnwidth]{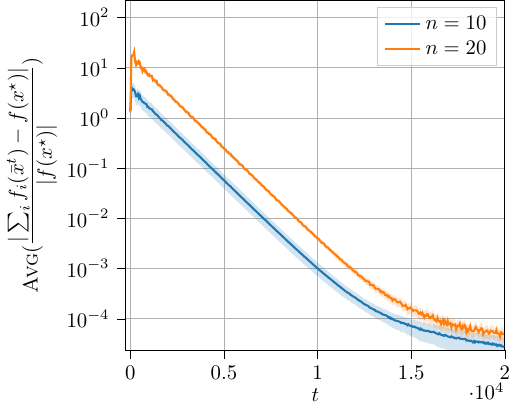}
	\includegraphics[width=.48\columnwidth]{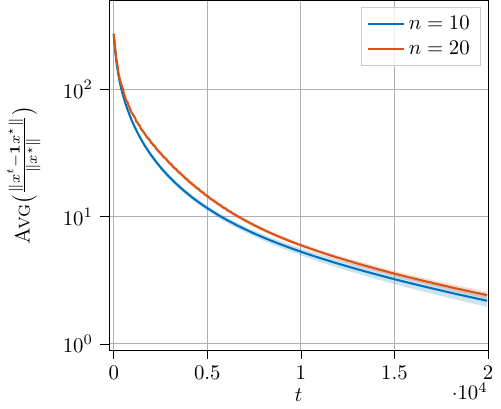}
	\caption{Monte Carlo simulations for large scale problems: mean and $1$-standard deviation band of cost error (left) and variable error (right).}\label{fig:cost_error_larger}
\end{figure}

First, we perform numerical simulations over large-scale problems with networks made of $N=250$
agents. 
We consider different optimization variable sizes, namely
$\xD=10,20$. 
Part of these simulations has been run on the Marconi100 HPC Cluster of the Italian Cineca. 
We used $10$ nodes of the cluster and, for each
node, we used $25$ cores and $4$ GPUs. 
The code has been adapted in order to
perform part of the computation directly on GPUs. 
The results are shown in Fig.~\ref{fig:cost_error_larger}. %
In detail, as the number of agents increases, the Lipschitz constant of the system to be averaged increases too.
Moreover, a larger domain of initial conditions also implies a potentially larger $L_\phi$ constant (cf.~\eqref{eq:L_phi}). %
This implies smaller $\gamma^\star$, which, fixed the other parameters, makes the convergence slower. 
The decision variable dimension instead impacts the selection of the dither signal. 
A larger number of states implies a larger number of frequencies. 
This, in turn, means a longer time to estimate the gradient (cf. Lemma~\ref{lemma:ESApprox}). 
Notice that, however, the accuracy of the final estimate is guaranteed by design. 
Indeed, since $\delta$ and $\gamma$ are designed on $\bar{\rho}$, the trajectories of (5) converge to a ball of radius $\bar{\rho}$ independently of the problem size.

  \subsection{Monte Carlo simulations varying number of agents}

\begin{figure}
    \centering
    \includegraphics[width=.4975\columnwidth]{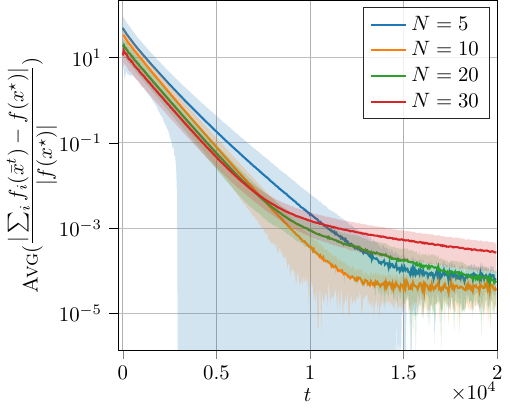}
    \includegraphics[width=.48\columnwidth]{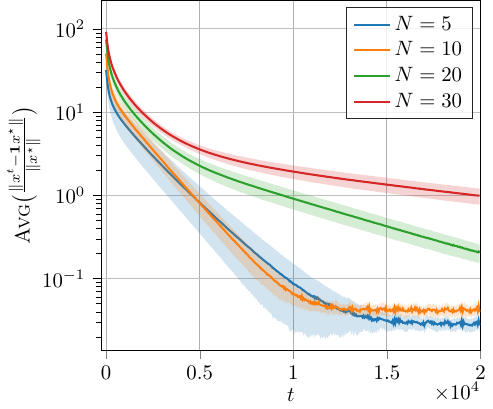}
    \caption{Monte Carlo simulations for a varying number of agents: mean and $1$-standard deviation band of cost error (left) and variable error (right).}\label{fig:cost_error}
  \end{figure}

  Second, we test \algo/ in a framework with $\xD = 10$ and different number of agents, i.e., $\nAg =5,10,20,30$.
  For each value of $N$, we generate communication graphs with a diameter $d$ such that the ratio $N/d$ is constant while varying $N$.
  The results depicted in Fig.~\ref{fig:cost_error} show that the algorithm slows down as the number of agents increases.

\subsection{Monte Carlo simulations varying the size $\xD$}

\begin{figure}[H]
    \centering
    \includegraphics[width=.4975\columnwidth]{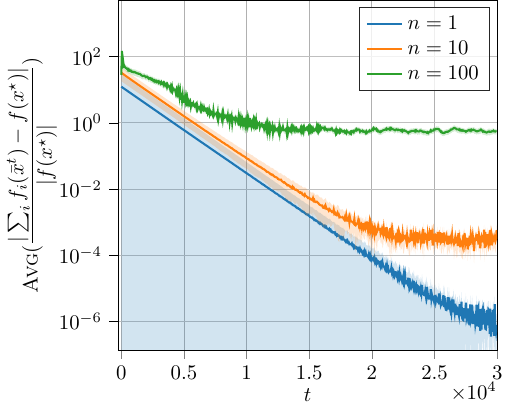}
    \includegraphics[width=.48\columnwidth]{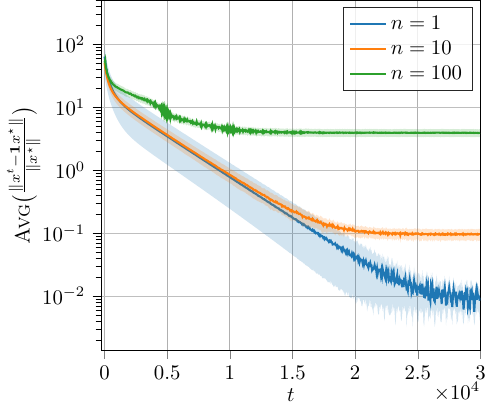}
    \caption{Monte Carlo simulations varying $\xD$: mean and $1$-standard deviation band of cost error (left) and variable error (right).}\label{fig:cost_error_1_10_100}
  \end{figure}

\begin{figure}[H]%
	  \centering
	  \includegraphics[width=.5\columnwidth]{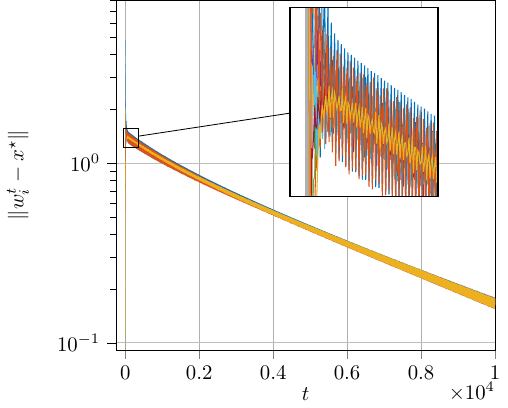}
	  \caption{Evolution of the agents' estimates $w_i^t$ in coordinate error with respect to the optimal solution $x^\star$.}
	  \label{fig:single_agent_trajectories}
  \end{figure}

Then, we test the algorithm features with a fixed number of agents $\nAg = 10$ and different decision variable dimensions $\xD = 1, 10, 100$.
The achieved results are reported in Fig.~\ref{fig:cost_error_1_10_100} and show that Algorithm~\ref{table:algorithm} slows down as $\xD$ increases.
  We conclude this part by providing in Fig.~\ref{fig:single_agent_trajectories} the results of a single trial in the case with $\xD = 1$ to show the evolution of the solution estimate of each agent in error coordinate with respect to the optimal solution to the problem.

\subsection{Monte Carlo simulations varying the parameter $\delta$}

\begin{figure}[H]%
    \centering
    \includegraphics[width=.4975\columnwidth]{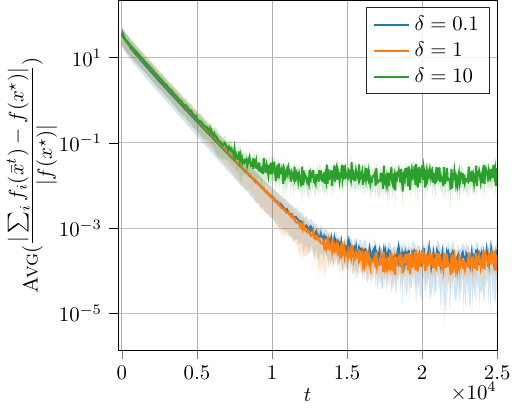}
    \includegraphics[width=.48\columnwidth]{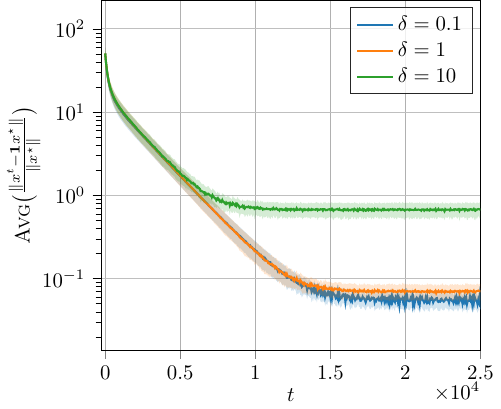}
    \caption{Monte Carlo simulations for varying amplitude $\delta$: mean and $1$-standard deviation band of cost error (left) and variable error (right).}\label{fig:monte_delta}
\end{figure}

Third, we provide numerical simulations in which we vary $\delta$ in the case in which $\nAg = \xD = 10$.
As one may expect from Theorem~\ref{th:MainResult}, we provide Fig.~\ref{fig:monte_delta} to show that the final accuracy of the algorithm increases with smaller values of $\delta$.

\subsection{Monte Carlo simulations for comparisons}
\begin{figure}[H]
	\centering
	\includegraphics[width=.4975\columnwidth]{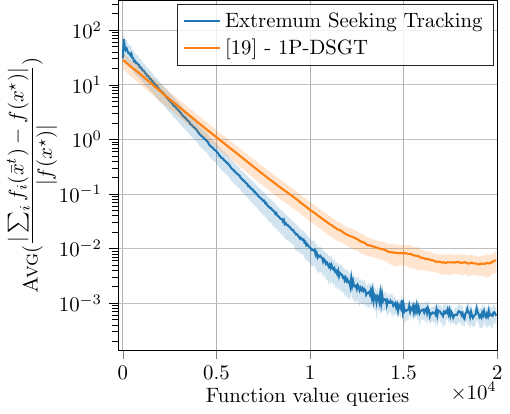}
	\includegraphics[width=.48\columnwidth]{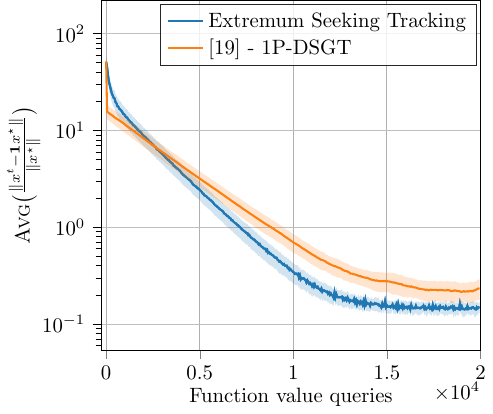}
	\caption{Monte Carlo simulations for comparison between \algo/ and 1P-DSGT in~\cite{mhanna2022zero}: mean and $1$-standard deviation band of cost error (left) and variable error (right).}\label{fig:comparison}
\end{figure}

We now perform simulations with $\xD=30$ and $\nAg=10$ to compare our method with the 1-Point Distributed Stochastic Gradient-Tracking Method (1P-DSGT) by~\cite{mhanna2022zero}.
We remark that, similarly to our scheme, when running 1P-DSGT each agent estimates the local gradient with one query of the objective function at each iteration.
These results are depicted in Fig.~\ref{fig:comparison} and have been obtained by running the considered distributed schemes with the same objective functions, communication graphs, and initial conditions of the solution estimates.

In detail, Fig.~\ref{fig:comparison} shows that, although the algorithm by~\cite{mhanna2022zero} exhibits a faster convergence in the beginning phase of the simulations, our scheme has a better convergence rate and final accuracy.
Both 1P-DSGT and \algo/ employ a single function query per agent to approximate the gradient. These methods rely on the concept that averaging across iterations around the \lq\lq quasi-static" local solution estimate $x_i^t$ provides an accurate gradient approximation. However, in the context of their theoretical frameworks, 1P-DSGT attains this approximation using an infinite number of samples (mean of an ergodic process), whereas \algo/ only needs $\tau_{\text{per}}$ samples.
This makes the convergence rate of our algorithm faster.

\subsection{Monte Carlo simulations in stochastic scenarios}

\begin{figure}%
    \centering
    \includegraphics[width=.4975\columnwidth]{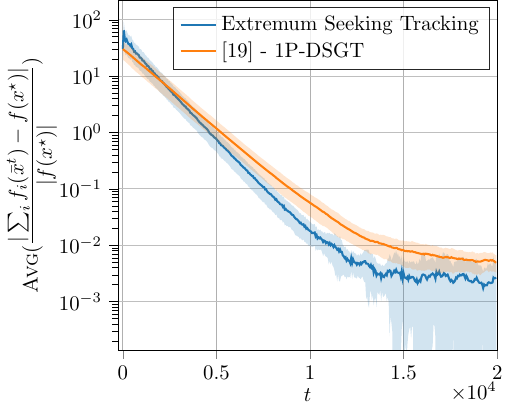}
    \includegraphics[width=.48\columnwidth]{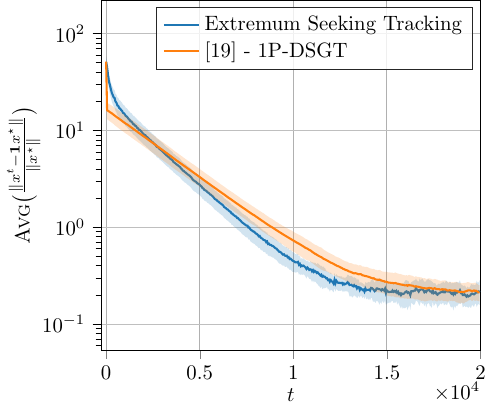}
    \caption{Monte Carlo simulations in stochastic scenarios: mean and $1$-standard deviation band of cost error (left) and variable error (right).}\label{fig:noisy}
  \end{figure}

Finally, we compare the considered distributed algorithms in a stochastic scenario in which each agent receives cost measurements affected by noise, i.e., in the case of $f_i(w_i^t) + \eta_i^t$ in place of $f_i(w_i^t)$ in Algorithm~\ref{table:algorithm}, where each component $\eta_i^t \in \real^{\xD}$ is randomly generated according to the Gaussian distribution with expected value $0$ and standard deviation $0.1$.
Fig.~\ref{fig:noisy} compares the behavior of \algo/ and 1P-DSGT by~\cite{mhanna2022zero} in the case in which $\nAg = 30$ and $\xD = 10$.
These plots confirm that \algo/ exhibits faster convergence and greater final accuracy with respect to the considered scheme also in the considered stochastic setup.
To interpret these results, we also remark that Theorem~\ref{th:MainResult} ensures (semi-global, practical) stability properties for \algo/.
Coherently, Fig.~\ref{fig:noisy} shows the typical behavior exhibited by the trajectories of perturbed systems in the neighborhood of (practically) stable equilibria.

\section{Conclusions}
\label{sec:Conclusion}

In this paper, we addressed a distributed optimization problem in which the cost
function is unknown and agents have only access to local measurements.  
Taking inspiration from a continuous gradient tracking algorithm, we proposed a
novel gradient-free distributed optimization algorithm in which gradients are
estimated via extremum seeking.
We analyzed the convergence properties of the proposed algorithm by using Lyapunov and averaging tools from system theory. 
We corroborated the theoretical analysis through Monte Carlo simulations on personalized optimization problems.

\appendices

\begin{appendix}

\subsection{Proof of Lemma~\ref{lemma:ESApprox}}
\label{app:ProofLemmaESApprox}

Given
    $\alpha =\col(\alpha_1,\dots,\alpha_\xD) \in \mathbb{N}^\xD$,
    $y = \col(y_1,\dots,y_\xD) \in\mathbb{R}^{\xD}$, and a smooth function
    $f\,:\,\mathbb{R}^\xD \to \mathbb{R}$, we define %
\begin{align*}
	\alpha ! &:= \alpha_1 !  \dots \alpha_\xD!,
	\quad
	&&y^\alpha := y_1^{\alpha_1}  \dots  y_\xD^{\alpha_\xD},
	\\
	\partial^\alpha f(y) &:= 
	\dfrac{\partial^{\alpha_1}}{\partial y_1^{\alpha_1}}	\dots\dfrac{\partial^{\alpha_\xD}}{\partial
                               y_n^{\alpha_\xD}} f(y),
	\quad
	&&|\alpha| := \alpha_1 + \dots + \alpha_\xD.
\end{align*}
Being each function $f_i$ smooth (cf. Assumption~\ref{ass:lipschitz}), we can apply Taylor's expansion (cf.~\cite[Theorem~2]{folland2005higher}) and write
\begin{align}
	f_i(x_i+\delta d_i^t) \!&= \!f_i(x_i) \! + \! \delta   {d_i^t}^\top \nabla f_i(x_i)  
	\! + \! \frac{\delta^2}{2}  {d_i^t}^\top  \nabla^2 f_i(x_i) d_i^t
	\notag\\
	&\hspace{.4cm}
	+ \delta^3 R_{i, 2}(x_i,\delta d_i^t),%
\label{eq:taylor_expansion}
\end{align}
where the remainder $R_{i, 2}(x_i, \delta d_i^t)$ is given by
\begin{align}
	R_{i, 2}(x_i,\delta d_i^t) = \sum_{|\alpha| = 3} \dfrac{\partial^\alpha f_i(x_i + c\delta d_i^t)}{\alpha!}(\delta d_i^t)^\alpha,\label{eq:remainder}
\end{align}
for some $c \in (0,1)$.
Then, we can use~\eqref{eq:taylor_expansion} to write
\begin{align}
	&\dfrac{2}{\delta\per}\sum_{k=t+1}^{t+\per}  d_i^k  f_i(x_i+\delta d_i^k)  
\notag
\\
	&=\dfrac{2 f_i(x_i)}{\delta\per} \sum_{k=t+1}^{t+\per} d_i^k  + \left[\dfrac{2}{\per} \sum_{k=t+1}^{t+\per}  \left(d_i^k  {d_i^k}^\top\right)\right] \nabla f_i(x_i)  
\notag\\
	&\hspace{.4cm}  + \dfrac{\delta}{\per} \sum_{k=t+1}^{t+\per}  \left(d_i^k {d_i^k}^\top\right) \nabla^2 f_i(x_i) d_i^k 
\notag\\
	&\hspace{.4cm}
	+  \dfrac{2}{\delta\per} \sum_{k=t+1}^{t+\per} d_i^k R_{i, 2}(x_i,\delta d_i^k).\label{eq:sum}
\end{align}
By combining~\eqref{eq:frequencies} and~\eqref{eq:sum}, we get
\begin{align*}
	&\sum_{k=t+1}^{t+\per}  d_i^k = 0, 
	\quad
	\dfrac{2}{\per}\sum_{k=t+1}^{t+\per}  \left(d_i^k  {d_i^k}^\top\right) = I_\xD
	\\
	&\sum_{k=t+1}^{t+\per}  \left(d_i^k  {d_i^k}^\top\right)\nabla^2 f_i(x_i) d_i^k = 0,
\end{align*}
which combined with~\eqref{eq:sum}, allow us to write 
\begin{align*}
		f_i(x_i+\delta d_i^t)= \nabla f_i(x_i) 
	+ \dfrac{2}{\delta\per} \sum_{k=t+1}^{t+\per} d_i^kR_{i, 2}(x_i,\delta d_i^k).
\end{align*}
The proof follows by setting $\ell_i(x_i) =\dfrac{2}{\per}\sum_{k=t+1}^{t+\per} d_i^k R_{i, 2}(x_i,\delta d_i^k)/\delta^3$. 
Finally, given a compact set $\mathcal{S}_i \subset \R^\xD$, let us bound $\norm{\ell_i(x_i)}$ for all $x_i \in \mathcal{S}_i$.
Note that $\norm{\delta d_i^t} \leq \delta\sqrt{\xD}$ for all $t \in \N$ and let $\mathcal{S}^\prime_i \subset \R^\xD$ be a compact set such that (i) $\mathcal{S}_i \subseteq \mathcal{S}_i^\prime \subset \R^\xD$, and (ii) $x_i + \delta d_i^t \in \mathcal{S}^\prime$ for all $x_i \in \mathcal{S}_i$, $\delta \in (0,1]$, and $t \in \N$.
Thus, we can write
\begin{align}
	&\sup_{\substack{x_i \in \mathcal{S}_i^\prime\\k \in \until{\per-1}}}\norm{\frac{R_{i, 2}(x_i,\delta d_i^k)}{\delta^3}} 
	\notag\\
	&\stackrel{(a)}{\leq} \sup_{\substack{x_i \in \mathcal{S}_i^\prime\\k \in \until{\per-1}}}\norm{\frac{1}{\delta^3}\sum_{|\alpha| = 3} \dfrac{\partial^\alpha f_i(x_i)}{\alpha!}(\delta d_i^k)^\alpha}
	\notag\\
	&\stackrel{(b)}{=}\sup_{\substack{x_i \in \mathcal{S}_i^\prime\\k \in \until{\per-1}}}\norm{\sum_{|\alpha| = 3} \dfrac{\partial^\alpha f_i(x_i)}{\alpha!}(d_i^k)^\alpha} 
	=: L_{i,\mathcal{S}_i}^\prime,\label{eq:L_i_S}
\end{align}
where in $(a)$ we use the expression~\eqref{eq:remainder} of $R_{i, 2}(x_i,\delta d_i^k)$, the definition of $\mathcal{S}_i^\prime$, and the fact that $\delta \in (0,1]$, while in $(b)$ we drop out the term $\delta^3$ from $(\delta d_i^k)^\alpha$.
We underline that, since the set $\mathcal{S}_i^\prime$ is compact and $f_i$ is smooth, $L_{i,\mathcal{S}_i}^\prime$ exists and is finite.
The bound of $\ell_i(x_i)$ follows by defining $L_{i,\mathcal{S}_i} := 2L_{i,\mathcal{S}_i}^\prime\sqrt{\xD}$ and combining the result~\eqref{eq:L_i_S} with the bound about the norm of the dither signal, i.e., $\norm{d_i^t} \leq \sqrt{\xD}$ for all $t \in \N$.

\subsection{Proof of Lemma \ref{lemma:AVunperturbed}}
\label{sec:ProofLemmaAVunperturbed}

		In~\cite{CARNEVALE2023110726}, it is provided a Lyapunov function proving that, under the Assumptions~\ref{ass:connectivity},~\ref{ass:strong_convexity}, and~\ref{ass:lipschitz}, the point $\xi^\star = (\1x^\star,\tilde{z}_\perp^{\text{eq}})$, with $\tilde{z}_\perp^{\text{eq}} := -R^\top G(\1 x^\star)$,  is a globally exponentially stable equilibrium for the continuous-time system 
		\begin{align*}
			\dot{\xi}(t) = \phi_{\text{GT}}(\xi(t)).
		\end{align*}
		In detail, \cite[Th.~3.1]{CARNEVALE2023110726} introduces a full-rank matrix $\bar{T} \in \R^{n_\xi \times n_\xi}$ to define $\bar{\xi} := \bar{T}\xi$, and a matrix $\bar{P} = \bar{P}^\top \in \R^{n_\xi \times n_\xi}$ such that 
		\begin{subequations}
			\begin{align}
				\bar{a}_1 I &\leq \bar{P} \leq \bar{a}_2 I\label{eq:P_positive_definite_bar}
				\\
				2\bar{\xi}\av^\top \bar{P} \bar{T}\phi_{\text{GT}}(\bar{T}^{-1}\bar{\xi}\av) 
				&\leq - \bar{a}_3\norm{\bar{\xi}\av}^2,\label{eq:cotninuous_negativeness_bar}
			\end{align}
		\end{subequations}
		for all $\xi\av \in \R^{n_{\xi}}$ and some $\bar{a}_1, \bar{a}_2, \bar{a}_3 > 0$.
		Let $P := \bar{T}^\top \bar{P}T$.
		Then, there exist $a_1, a_2, a_3 > 0$ such that
		\begin{subequations}
			\begin{align}
				a_1 I &\leq P \leq a_2 I\label{eq:P_positive_definite}
				\\
				2\xi\av^\top P \phi_{\text{GT}}(\xi\av)  &\leq - a_3\norm{\xi\av}^2,\label{eq:cotninuous_negativeness}
			\end{align}
		\end{subequations}
		for all $\xi\av \in \R^{n_{\xi}}$.
		Hence, we define the candidate Lyapunov function $V(\xi\av) := \xi\av^\top P \xi\av$ and bound $\Delta V(\xi\av^t) := V(\xi\av^{t+1}) - V(\xi\av^t)$ along the trajectories of~\eqref{eq:GTAUnperturbed} as
		\begin{align}
			\Delta V(\xi\av^t) \leq -\gamma a_3\norm{\xi\av^t}^2 + \gamma^2 \phi_{\text{GT}}(\xi\av^t)^\top P \phi_{\text{GT}}(\xi\av^t).\label{eq:deltaV_negative}
		\end{align}
		Moreover, by using the Lipschitz continuity of the gradients of the objective functions (cf. Assumption~\ref{ass:lipschitz}) and the definition of $\phi_{\text{GT}}$, there exists $a_4 > 0$ such that 
		\begin{align}
			\norm{\phi_{\text{GT}}(\xi\av^t)} \leq a_4\norm{\xi\av^t}.\label{eq:phi_GT_bound}
		\end{align}
	Finally, for any $c_1 \in (0,a_3)$, let $\gamma_0 := (a_3-c_1)/(a_2 a_4^2)$ and the proof follows by using~\eqref{eq:deltaV_negative} and~\eqref{eq:phi_GT_bound}.

\subsection{Proof of Lemma \ref{lemma:Average}}
\label{sec:ProofLemmaAverage}

	The proof relies on (i) the matrix $P$ satisfying~\eqref{eq:V_AVunperturbed}, and (ii) the fact that the norm of the perturbation term $\gamma \delta B u(\xi\av^t)$ can be arbitrarily reduced through the parameter $\delta$ as long as $\xi\av^t$ lies into a compact set.
	First of all, without loss of generality, we assume $\rho \leq r_\xi$. 
	Indeed, we will use the parameter $r_\xi$ to define a (compact) ball and arbitrarily bound the norm of the perturbation term $\gamma \delta B u(\xi\av^t)$ through the parameters $\delta$ as long as $\xi\av^t$ lies into this ball.
	Hence, we can always use the more conservative condition.
	In detail, we introduce the candidate Lypaunov function 
		$V(\xi\av) := \xi\av^\top P \xi\av$
	and the set 
	$\Omega_{r_{\xi}} := \{ \xi\av\in\real^{n_\xi}  \mid V(\xi\av) \leq a_2 r_{\xi}^2\}\subset \R^{n_\xi}$.
	Then, from~\eqref{eq:V_bounds}, we derive $\mathcal{B}_{r_{\xi}} \subseteq \Omega_{r_{\xi}} \subseteq \mathcal{B}_{r^\prime_{\xi\av}}$, where $r^\prime_{\xi\av} := \sqrt{a_2/a_1}r_{\xi}$.
	Now, under the assumption $\xi\av^t \in \mathcal{B}_{r_{\xi}}$ (later verified by a proper selection of the algorithm parameters), it holds $\xi\av^t \in \Omega_{r_{\xi}}$.
	By using this property and since $\gamma \in (0,\gamma_0)$, we use~\eqref{eq:V_negative}, the Cauchy-Schwarz inequality, and~\eqref{eq:phi_GT_bound} to bound $\Delta V(\xi\av^{t}) := V(\xi\av^{t+1}) - V(\xi\av^t)$ along the trajectories of~\eqref{eq:xia} as
	\begin{align}
	\Delta V(\xi\av^t) &\le  - \gamma c_1 \|\xi\av^t\|^2 +  \gamma \delta^2 2\norm{PB} \|\xi\av^t\|\norm{u(\xi\av^t)}
	\notag\\
	&\hspace{.4cm}
	+ \delta^2 \gamma^2 2a_4\norm{PB}\norm{\xi\av^t}\norm{u(\xi\av^t)}
	\notag\\
	&\hspace{.4cm}
	+ \delta^4 \gamma^2 \norm{B^\top PB}\norm{u(\xi\av^t)}^2.\label{eq:DeltaV_perturbed}
	\end{align}
	Now, for all $i \in \until{N}$, let us introduce the $\mathcal{S}_i := \{ x_i \in \R^\xD \mid \norm{x_i - x^\star} \leq r^\prime_{\xi\av}\}$.
	Hence, we note that $\xi\av := (\tilde{x}\av, \tilde{z}_{\perp,\text{avg}}) \in \Omega_{r_{\xi}} \implies  \tilde{x} \in \mathcal{S} \subset \R^{\nAg\xD}$, where $\mathcal{S} := \mathcal{S}_1 \times \dots \times \mathcal{S}_N$.
	Then, we apply result~\eqref{eq:ell_bound} to claim that, for all $i \in \until{N}$, it holds $\ell_i(x_i) \leq L_{i,\mathcal{S}_i}$ for all $x_i \in \mathcal{S}_i$.
	Thus, by defining $L_{\mathcal{S}} := \max_{i}\{L_{1,\mathcal{S}_1}, \dots, L_{N,\mathcal{S}_N}\}$ and using the definition $u(\xi\av^t) = \ell(\tilde{x}^t + \1x^\star)$, we get
	\begin{align}
		\norm{u(\xi\av^t)} \leq \sqrt{N}L_{\mathcal{S}}.\label{eq:u_bound}
	\end{align}
	Since $\gamma \in (0,1]$ and $\delta \in (0,1]$, we bound~\eqref{eq:DeltaV_perturbed} as
	\begin{align}
		\Delta V(\xi\av^t) &\le - \gamma c_1 \|\xi\av^t\|^2 +  \gamma \delta^2 (b_1\norm{\xi\av^t} + b_2),\label{eq:DeltaV_perturbed_2}
	\end{align}
	where we introduced
	\begin{align*}
		b_1 &:=  2\norm{PB} \sqrt{N}L_{\Omega_{r_{\xi}}} + 2a_4\norm{PB}\sqrt{N}L_{\mathcal{S}}%
		\\
		b_2 &:= \norm{B^\top PB} N L_{\mathcal{S}}^2.%
	\end{align*}
	Therefore, for any $\rho \in (0,r_\xi)$ and $c_1^\prime \in (0,c_1)$, we define 
	\begin{align}
		\delta^\star_1 := \min\left\{\sqrt{(c_1 - c_1^\prime)\rho^2/(b_1 r^\prime_{\xi\av} + b_2)},1\right\}.\label{eq:delta_star}
	\end{align}
	Hence, by combining~\eqref{eq:DeltaV_perturbed_2} and~\eqref{eq:delta_star}, if $\delta \in (0,\delta^\star_1)$, then, for all $\xi\av^t \in \Omega_{r^\prime_{\xi\av}}$ such that $\norm{\xi\av^t} \ge \rho$, it holds
	\begin{align}
		\Delta V(\xi\av^t) \leq -\gamma c_1^\prime \norm{\xi\av^t}^2.\label{eq:DeltaV_perturbed_3}
	\end{align}
	Thus, the inequality~\eqref{eq:DeltaV_perturbed_3} ensures that the set $\Omega_{r_{\xi}}$ is forward-invariant for system~\eqref{eq:xia}.
	Hence, if we pick $\xi\av^0 \in \mathcal{B}_{r_{\xi}}$, we prove that $\xi\av^{t} \in \Omega_{r_{\xi}}$ for all $t \in \N$.
	Consequently, the bound~\eqref{eq:u_bound} holds for all $t \in \N$ and, in turn, also the inequality~\eqref{eq:DeltaV_perturbed_3} is verified for all $t \in \N$, namely we proved that the trajectories of system~\eqref{eq:xia} enter the ball $\mathcal{B}_\rho$ exponentially fast.
	The result~\eqref{eq:exponential_lemma} follows from the inequality~\eqref{eq:DeltaV_perturbed_3} and~\eqref{eq:V_bounds} by setting $c_3 := c_1^\prime/(2a_2)$.

	\end{appendix}

\end{document}